\documentstyle[12pt]{article}

\newfont{\frak}{eufm10 scaled\magstep1}
\newfont{\extra}{msbm10 scaled\magstep1}

\newcommand{\sect}[1]{\setcounter{equation}{0}\section{#1}}

\begin{document}

\begin{center}
{\Large {\bf Quantum Jordanian twist }}\\[5mm] 
{Vladimir D. Lyakhovsky $^{a,b}$}\\ 
{Alexander M. Mirolubov $^b$} 
\\ {Mariano A. del Olmo $^a$}\\[5mm] 
{\sl $^a$ Departamento de F\'{\i}sica Te\'orica, Facultad de Ciencias\\
Universidad de Valladolid, E-47011, Valladolid, Spain. \\[2mm] 
$^b$ Theoretical Department, St. Petersburg State University,\\ 
198904, St. Petersburg, Russia}
\end{center}

\vskip 0.5cm
%%%%%%%%%%%%%%%%%%%%%%%%%%%%%%%%%%%%%%%%%%%%%%%%%%%%%%%%%%%%%%%%%%%%%%%
\begin{abstract}
The quantum deformation of the Jordanian twist ${\cal F}_{q{\cal J}}$ for
the standard quantum Borel algebra $U_q (B)$ is constructed. It gives the
family $U_{q {\cal J}}(B)$ of quantum algebras depending on parameters $\xi$
and $h$. In a generic point these algebras represent the hybrid
(standard--nonstandard) quantization. The quantum Jordanian twist can be
applied to the standard quantization of any Kac--Moody algebra. The
corresponding classical $r$--matrix is a linear combination of the
Drinfeld--Jimbo and the Jordanian ones. The two-parametric families
of Hopf algebras  obtained here are smooth and for the limit values of the parameters the
standard and nonstandard quantizations are recovered. The twisting element $%
{\cal F}_{q{\cal J}}$ also has the correlated limits, in particular when $q$
tends to  unity it acquires the canonical form of the Jordanian twist. To
illustrate the properties of the quantum Jordanian twist we construct the
hybrid quantizations for $U (sl(2))$ and for the corresponding affine
algebra $U (\widehat{ sl(2)})$. The universal quantum ${\cal R}$--matrix and
its defining representation are presented.
\end{abstract}
%%%%%%%%%%%%%%%%%%%%%%%%%%%%%%%%%%%%%%%%%%%%%%%%%%%%%%%%%%%%%%%%%%%%%%
%%%%%%%%%%%%%%%%%%%%%%%%%%%%%%%%%%%%%%%%%%%%%%%%%%%%%%%%%%%%%%%%%%%%%%

\sect{Introduction}

It is known for a long time \cite{DRI} that a Hopf algebra ${\cal A}%
(m,\Delta,\epsilon,S)$ with multiplication $m\colon {\cal A}\otimes {\cal A}%
\to {\cal A}$, coproduct $\Delta \colon {\cal A}\to {\cal A} \otimes {\cal A}
$, counit $\epsilon \colon {\cal A}\to C$, and antipode $S : {\cal A}\to
{\cal A}$ can be transformed with an invertible (twisting) element ${\cal F}%
\in {\cal A} \otimes {\cal A}$, ${\cal F}=\sum f_i^{(1)}\otimes f_i^{(2)}$,
into a twisted one ${\cal A}_{{\cal F}}(m,\Delta _{{\cal F}},\epsilon ,S_{%
{\cal F}})$ that have the same multiplication and counit but different
coproduct and antipode. The twisted coproduct is given by
\begin{equation}
\label{def-t} \Delta _{{\cal F}}(a)={\cal F}\Delta (a){\cal F}^{-1}.
\end{equation}
The twisting element has to satisfy the equations
\begin{eqnarray}
\label{def-n}
(\epsilon \otimes  id)({\cal F}) = (id \otimes  \epsilon)({\cal F})=1,
\\[0.2cm]
\label{gentwist}
{\cal F}_{12}(\Delta \otimes  id)({\cal F}) =
{\cal F}_{23}(id \otimes  \Delta)({\cal F}).
\label{TE}
\end{eqnarray}
There are several special types of twists. For our purposes the most
interesting will be the factorizable twist whose twisting element satisfies
the factorized twist equations \cite{RES}:
\begin{equation}
\label{fact}
\begin{array}{l}
(\Delta \otimes id)(
{\cal F}) = {\cal F}_{13}{\cal F}_{23}, \\ [0.2cm] ( id \otimes \Delta_{%
{\cal F}})({\cal F}) = {\cal F}_{12}{\cal F}_{13}.
\end{array}
\end{equation}
If the initial Hopf algebra ${\cal A}$ is quasitriangular with universal $%
{\cal R}$--matrix ${\cal R}$ then ${\cal A}_{{\cal F}}$ is the twisted Hopf
algebra whose universal element ${\cal R}_{{\cal F}}$ is related to the
initial one by
\begin{eqnarray}
\label{RF}
{\cal R}_{\cal F}=  {\cal F}_{21} \,{\cal R} \,{\cal F}^{-1}.
\end{eqnarray}

The Jordanian twist with the two-dimensional carrier subalgebra $B(2)$,
$$
\left[ H,E \right] =E,
$$
defined by the canonical twisting element
\begin{equation}
\label{j-twist} {\cal F}_{{\cal J}}^{\ c} = e^{H \otimes \sigma}, \qquad
\sigma = \ln (E + 1),
\end{equation}
is the first and very important example \cite{OGI} of a nontrivial twist
with explicitly defined twisting element.

It was proved in \cite{OGI} that there exist mixed quantizations combining
the properties of the standard deformations and that of the twisted
algebras. Up to now a considerable amount of studies devoted to combined
(standard--nonstandard) quantizations were performed (especially for the
case of $U(sl(2))$, $U(gl(2))$ and the corresponding quantum groups). See,
for example, the works by Gerstenhaber {\sl et al} \cite{GGS}, Kupershmidt
\cite{KUP}, Ballesteros {\sl et al} \cite{BHP}, Abdesselam {\sl et al} \cite
{CHA,ABD}, Aneva {\sl et al} \cite{CHA2} and references therein. (The last
work contains a kind of review of the situation with the combined
deformations and we shall return to it in Section \ref{conclusions}.) But a
question remains, whether it is possible to supply the combined quantization
with a twisting element that would bring it ``back'' to the standard quantum
algebra (such as $U_q(sl(2))$). This ``going back'' procedure is a limit
process. In fact there are two limits to be considered, they are related
with the behaviour of two main parameters: the deformation parameter $h=\ln q$
and the twisting parameter $\xi $. Recently a $q$--analog $({\cal F}_{{\cal J}%
}^{\ c})_q$ of the Jordanian twisting element (\ref{j-twist}) was
constructed \cite{KTS}. It transforms the standard quantization $U_q(sl(2))$
into the combined deformation $(U_q(sl(2)))_{{\cal J}}$ and the inverse
operator $({\cal F}_{{\cal J}}^{\ c})_q^{-1}$ obviously brings the algebra $%
(U_q(sl(2)))_{{\cal J}}$ back into the standard deformation. For us it is
important to notice that the Hopf algebra $(U_q(sl(2)))_{{\cal J}}$ has no
classical limit for $h\to 0$.

In this paper we demonstrate that there are other sheets of combined
quantizations for which both limits exist: for $h \to 0$ (the standard $q$%
--deformation) and for $\xi \to 0$ (the nonstandard or Jordanian
deformation). We investigate the existence of quantum deformations that do
not only refer to the combined classical $r$--matrix and can be connected
with the standard quantization by a twist, but such that all their algebraic
elements (bialgebraic structure, universal ${\cal R}$--matrix and twisting
element) have well defined limits. In Section \ref{jordaniantwist} we
demonstrate that this problem can be solved by constructing a quantum
deformation ${\cal F}_{q{\cal J}}$ of the Jordanian twist ${\cal F}_{{\cal J}%
}$. This new quantum Jordanian twist ${\cal F}_{q {\cal J}}$ acts on the
standard quantizations of the universal enveloping algebras and transforms
them into the hybrid quantizations (we have borrowed this term from Ref.
\cite{CHA2}). The twisting element ${\cal F}_{q {\cal J}}$ itself and all
the corresponding twisted constructions have both natural limits. In Section
\ref{jordanian deformations} we apply this twist to the quantum algebras
based on $sl(2)$ and get the hybrid quantizations of $U(sl(2))$ and of the
quantum affine algebra $U(\widehat{sl(2)})$,
$$
\begin{array}{l}
U_q(sl(2))
\stackrel{{\cal F}_{q {\cal J}}}{\longrightarrow} U_{q{\cal J}}(sl(2)), \\
[0.2cm] U_q(\widehat{sl(2)}) \stackrel{{\cal F}_{q {\cal J}}}{\longrightarrow%
} U_{q{\cal J}}(\widehat{sl(2)}) .
\end{array}
$$
The corresponding universal ${\cal R}$-matrices and their defining
representations are also presented.
%%%%%%%%%%%%%%%%%%%%%%%%%%%%%%%%%%%%%%%%%%%%%%%%%%%%%%%%%%%%%%%%%%%%%%%%%%%%
%%%%%%%%%%%%%%%%%%%%%%%%%%%%%%%%%%%%%%%%%%%%%%%%%%%%%%%%%%%%%%%%%%%%%%%%%%%%

\sect{Quantum Jordanian twist for $U_q(B)$}

\label{jordaniantwist}

\newtheorem{Proposition}{Proposition}
\begin{Proposition}\label{proposition1}
The quantum Borel algebra $U_q (B)$
\begin{equation}
\label{q-bor}
\left[ H,E \right] =E, \qquad
\begin{array}{l}
\Delta (H) = H \otimes 1 + 1 \otimes H,\\[0.2cm]
\Delta (E) = E \otimes 1 + e^{hH} \otimes E,
\end{array}
\end{equation}
admits the twist with the element
\begin{equation}
\label{q-twist}
\widetilde{\cal F}_{q{\cal J}}  = e^{H \otimes \sigma},
\qquad \sigma = \ln (E + e^{hH}).
\end{equation}
\end{Proposition}
%%%%
{\bf Proof.} \thinspace  We shall demonstrate that the element (\ref{q-twist}%
) satisfies the factorized twist equations. The first of them is trivially
fulfilled due to the primitivity of $H$. To check the second let us change
the basis. The new generator
\begin{equation}
\label{brev-e}\breve E=E-1+e^{hH},
\end{equation}
has the same coproduct as $E$. Performing the substitution we get for $U_q(B)
$ the relations:
\begin{equation}
\label{brevq-bor}\left[ H,\breve E\right] =\breve E+1-e^{hH},\qquad
\begin{array}{l}
\Delta (H)=H\otimes 1+1\otimes H, \\
[0.2cm]\Delta (\breve E)=\breve E\otimes 1+e^{hH}\otimes \breve E,
\end{array}
\end{equation}
and $\sigma $ from (\ref{q-twist}) gets the form
$$
\sigma =\ln (1+\breve E).
$$
The adjoint action of $H$ on $\breve E$ differs from that of $H$ on $E$ by
terms that are central. So, one gets
$$
e^{{\rm ad}(H\otimes \sigma )}\circ \breve E\otimes 1=\breve E\otimes
e^\sigma +\left( 1-e^{hH}\right) \otimes \left( e^\sigma -1\right) .
$$
Now we can obtain the final form of the coproduct for $e^\sigma $,
\begin{equation}
\begin{array}{lcl}
\widetilde{\Delta }_{q{\cal J}}(e^\sigma ) & = 
& \widetilde{\Delta}_{q
{\cal J}}\left( \breve E+1\right)  \\ [0.2cm] & = & e^{
{\rm ad}(H\otimes \sigma )}\circ \left( \breve E\otimes 1+e^{hH}\otimes
\breve E+1\otimes 1\right)  \\ [0.2cm] & = & e^{
{\rm ad}(H\otimes \sigma )}\circ \left( \breve E\otimes 1\right)
+e^{hH}\otimes \breve E+1\otimes 1 \\ [0.2cm] & = & \breve E\otimes e^\sigma
+1\otimes e^\sigma = e^\sigma \otimes e^\sigma .
\end{array}
\end{equation}
Consequently, $\sigma $ becomes primitive with respect to the deformed
coproduct $\widetilde{\Delta}_{q{\cal J}}$. 
Thus, the element $\widetilde{\cal F}_{q%
{\cal J}}$ satisfies also the second of the factorized twist equations (\ref
{fact}). This completes the proof. $\spadesuit $ 
%%%%%%%%%%%%%%%%%%%%%%%%%

It will be useful to introduce now the twisting parameter $\xi$. This can be
achieved by rescaling the generator $E \to \xi E $. So, Proposition 1 is
valid also for the parametric set of twisting elements
\begin{equation}
\label{qxi-twist} \widetilde{\cal F}_{q{\cal J}} (h,\xi) = e^{H \otimes
\sigma}, \qquad \sigma = \ln (\xi E + e^{hH}).
\end{equation}
%%%%%%%%%%%%%%%%%%%%%%%%%

Performing the twisting we obtain the smooth two-parameter set $\{\widetilde{%
U}_{q{\cal J}}(B)(h,\xi )\}$ of Hopf algebras
\begin{equation}
\label{qjt-bor}\left[ H,E\right] =E,\qquad
\begin{array}{l}
\widetilde{\Delta }_{q{\cal J}}(H)=H\otimes 1+e^{{\rm ad}(H\otimes \sigma
)}\circ (1\otimes H), \\ [0.2cm]\widetilde{\Delta }_{q{\cal J}}(E)=\frac
1\xi \left( e^\sigma \otimes e^\sigma -e^{{\rm ad}(H\otimes \sigma )}\circ
\left( e^{hH}\otimes e^{hH}\right) \right) .
\end{array}
\end{equation}
This set describes the {\sl hybrid quantization} that have the properties
both of the standard quantization and of the Jordanian. For the limit values
of the parameters we get the Hopf algebras whose characteristics are
different from the generic ones. The boundary of the set corresponding to $%
h=0$ gives the ordinary Jordanian quantization $\widetilde{U}_{q{\cal J}%
}(B)(0,\xi )=U_{{\cal J}}(B)$. The boundary $\widetilde{U}_{q{\cal J}%
}(B)(h,0)$ describes the Hopf algebras that are equivalent to the standard
deformation $U_q(B)$ but have the shifted coproduct for $E$. Such behaviour
is in accordance with the limit form of the twisting element (\ref{qxi-twist}%
)
\begin{equation}
\label{r-twist}\widetilde{\cal F}_{q{\cal J}}(h,0)=e^{hH\otimes H}.
\end{equation}

If we want to go back (when $\xi$ tends to zero) to the initial algebra $%
U_{q} (B)$ the compensating Reshetikhin twist must be applied. The final
construction is defined as follows. 
%%%%%%%%%%%%%%%%%%%%%%%%%

\begin{Proposition}\label{proposition2}
The quantum Borel algebra $U_q (B)$
$$
\left[ H,E \right] =E, \qquad
\begin{array}{l}
\Delta (H) = H \otimes 1 + 1 \otimes H,\\[0.2cm]
\Delta (E) = E \otimes 1 + e^{hH} \otimes E,
\end{array}
$$
admits the twist with the element
\begin{equation}
\label{qr-twist}
{\cal F}_{q{\cal J}}(h,\xi) = e^{H \otimes \omega} e^{-hH \otimes H},
\qquad \omega = \ln \left( \left( \xi E + 1 \right) e^{hH} \right).
\end{equation}
\end{Proposition}
%%%%%%%%%%%%%%%%%%%%%%%%%
{\bf Proof.} \, The first factor of the twisting element (\ref{qr-twist})
transforms the quantized Borel $U_q (B)$ into its dual Hopf algebra:
\begin{equation}
\label{qr-bor} \left[ H,E \right] =E, \qquad
\begin{array}{l}
\Delta_{qR} (H) = H \otimes 1 + 1 \otimes H, \\
[0.2cm] \Delta_{qR} (E) = E \otimes e^{-hH} + 1 \otimes E.
\end{array}
\end{equation}
The substitution $E \to E e^{-hH} $ brings us back to the initial Borel in
the form (\ref{q-bor}) and also changes $\omega$ in (\ref{qr-twist}) for $%
\sigma$ (as defined in (\ref{qxi-twist})). After this we find ourselves in
the situation of the Proposition \ref{proposition1}. $\spadesuit$
%%%%%%%%%%%%%%%%%%%%%%%%%

Applying the twist (\ref{qr-twist}) to $U_q (B)$ we get the Hopf algebras
with the following defining relations:
\begin{equation}
\label{qj-bor}
\begin{array}{l}
\qquad \qquad\qquad \left[ H,E \right] =E, \\
[0.2cm] \Delta_{q
{\cal J}} (H) = H \otimes 1 + e^{{\rm ad}(H \otimes \omega)} \circ (1
\otimes H), \\ [0.2cm] \Delta_{q{\cal J}} (E) = \frac{1}{\xi} \left( \left(
e^{\omega} \otimes e^{\omega} \right) - e^{{\rm ad}(H \otimes \omega)} \circ
\left( e^{-hH} \otimes e^{-hH} \right) - 1 \otimes 1 \right).
\end{array}
\end{equation}

We have two correlated smooth sets: the set $\left\{ U_{q{\cal J}}(B)(h,\xi
)\right\} $ of hybrid quantizations (\ref{qj-bor}) and the set of quantum
Jordanian twists $\left\{ {\cal F}_{q{\cal J}}(h,\xi )\right\} $. For each
point $U_{q{\cal J}}(B)(h,\xi )$ of the first set there exits a twist ${\cal %
F}_{q{\cal J}}^{-1}(h,\xi )$ that connects this point with the algebra $U_{q%
{\cal J}}(B)(h,0)$. Now the sets have the appropriate boundary behaviour:
\begin{equation}
\begin{array}{ccccc}
&  & U_{q{\cal J}}(B)(h,\xi ) &  &  \\
& \rule{25mm}{0mm}{\scriptsize h\rightarrow 0}\swarrow  &  & \searrow
{\scriptsize \xi \rightarrow 0}\rule{25mm}{0mm} &  \\
& U_{{\cal J}}(B)=U_{q{\cal J}}(B)(0,\xi ) &  & U_{q{\cal J}}(B)(h,0)=U_q(B)
&
\end{array}
\end{equation}
\smallskip

\begin{equation}
\begin{array}{ccccc}
&  & {\cal F}_{q{\cal J}}(h,\xi) &  &  \\
& \rule{25mm}{0mm} {\scriptsize h \rightarrow 0} \swarrow &  & \searrow
{\scriptsize \xi \rightarrow 0}\rule{25mm}{0mm} &  \\
& {\cal F}_{{\cal J}}^{\ c} = {\cal F}_{q{\cal J}} (0,\xi) &  & {\cal F}_{q%
{\cal J}}(h,0) = 1 \otimes 1 &
\end{array}
\end{equation}

For the internal points of $\left\{ U_{q{\cal J}} (B)(h,\xi) \right\} $ the
defining relations (\ref{qj-bor}) can be written in a compact form
\begin{equation}
\label{qjsi-bor}
\begin{array}{l}
\qquad \qquad \left[ H,e^{\omega} \right] =e^{\omega} -e^{hH}, \\
[0.2cm] \Delta_{q
{\cal J}} (H) = H \otimes 1 + e^{{\rm ad}(H \otimes \omega)} \circ (1
\otimes H), \\ [0.2cm] \Delta_{q{\cal J}} (e^{\omega}) = e^{\omega} \otimes
e^{\omega}.
\end{array}
\end{equation}
But such description becomes incomplete on the boundary $\left\{ U_{q{\cal J}%
} (B)(h,0) \right\} $ where $\omega(h,0) = hH$.
%%%%%%%%%%%%%%%%%%%%%%%%%%%%%%%%%%%%%%%%%%%%%%%%%%%%%%%%%%%%%%%%%%%%%%
%%%%%%%%%%%%%%%%%%%%%%%%%%%%%%%%%%%%%%%%%%%%%%%%%%%%%%%%%%%%%%%%%%%%%%%

\sect{Jordanian deformations of quantum algebras}

\label{jordanian deformations} %%%%%%%%%%%%%%%%%%%%%%%%%%
%%%%%%%%%%%%%%%%%%%%%%%%%%

\subsection{Hybrid quantization $U_{q{\cal J}} (sl(2))$}

The quantum Jordanian twists ${\cal F}_{q{\cal J}}(h,\xi)$ and $\widetilde{%
\cal F}_{q{\cal J}}(h,\xi) $ can be applied to any Hopf algebra containing
the quantized Borel algebra (\ref{q-bor}). Let us start with the standard
quantization $U_q(sl(2))$:
\begin{equation}
\label{q-sl2}
\begin{array}{l}
\left[ H,E_{\pm} \right] =\pm E_{\pm}, \\
[0.25cm] \left[ E_+,E_- \right] = \frac{e^{hH} - e^{hH}}{1 - e^h},
\end{array}
\qquad
\begin{array}{l}
\Delta_q (H) = H \otimes 1 + 1 \otimes H, \\
[0.2cm] \Delta_q (E_+) = E_+ \otimes 1 + e^{hH} \otimes E_+, \\
[0.2cm] \Delta_q (E_-) = E_- \otimes e^{-hH} + 1 \otimes E_-. \\
[0.2cm]
\end{array}
\end{equation}
Applying the twist ${\cal F}_{q{\cal J}}(h,\xi)$ in the form given in  (\ref{qr-twist}%
) we get the two-parameter set $\left\{ U_{q{\cal J}} (sl(2))(h,\xi)
\right\} $ of quantum deformations that are the {\sl hybrids} of standard
and Jordanian ones:
\begin{equation}
\label{qj-sl2}
\begin{array}{l}
\left[ H,E_{\pm} \right] =\pm E_{\pm}, \\
[0.2cm] \left[ E_+,E_- \right] =
\frac{e^{hH} - e^{hH}}{1 - e^h}, \\ [0.2cm] \Delta_{q
{\cal J}} (H) = H \otimes 1 + e^{{\rm ad}(H \otimes \omega)} \circ (1
\otimes H), \\ [0.2cm] \Delta_{q
{\cal J}} (E_+) = \frac{1}{\xi} \left( e^{\omega} \otimes e^{\omega} - e^{%
{\rm ad}(H \otimes \omega)} \circ \left( e^{-hH} \otimes e^{-hH} \right) - 1
\otimes 1 \right), \\ [0.2cm] \Delta_{q
{\cal J}} (E_-) = E_- \otimes e^{-\omega} + e^{{\rm ad}(H \otimes \omega)}
\circ \left( 1 \otimes E_- \right). \\ [0.2cm]
\end{array}
\end{equation}
The Hopf algebras $U_q(sl(2))$ and $U_{{\cal J}}(sl(2))$ form the boundaries
of this set:
\begin{equation}
\begin{array}{ccccc}
&  & U_{q{\cal J}} (B)(h,\xi) &  &  \\
& \rule{35mm}{0mm} {\scriptsize h \rightarrow 0} \swarrow &  & \searrow
{\scriptsize \xi \rightarrow 0}\rule{35mm}{0mm} &  \\
& U_{{\cal J}} (sl(2)) = U_{q{\cal J}} (sl(2))(0,\xi) &  & U_{q{\cal J}}
(sl(2))(h,0) = U_{q} (sl(2)) &
\end{array}
\end{equation}
When the internal subset $\left\{ U_{q{\cal J}} (sl(2))(h,\xi) \, | \, h>0,
\xi >0 \right\} $ is considered the compact form of the defining relations
can be used,
\begin{equation}
\label{qjsi-sl}
\begin{array}{l}
\left[ H,e^{\omega} \right] =e^{\omega} -e^{hH}, \\
[0.2cm] \left[ H, E_- \right] = -E_-, \\
[0.2cm] \left[ E_-,e^{\omega} \right]_{e^h} =\xi
\frac{1-e^{2hH}}{1-e^{-h}}, \\ [0.2cm] \Delta_{q
{\cal J}} (H) = H \otimes 1 + e^{{\rm ad}(H \otimes \omega)} \circ (1
\otimes H), \\ [0.2cm] \Delta_{q
{\cal J}} (e^{\omega}) = e^{\omega} \otimes e^{\omega}, \\ [0.2cm] \Delta_{q%
{\cal J}} (E_-) = E_- \otimes e^{-\omega} + e^{{\rm ad}(H \otimes \omega)}
\circ \left( 1 \otimes E_- \right).
\end{array}
\end{equation}

The algebra (\ref{q-sl2}) is quasitriangular with the universal ${\cal R}$%
-matrix
\begin{equation}
\label{uni-r} {\cal R}_q = e^{hH \otimes H} \sum_{n = 0}^{\inf} \frac{
\left(1- e^{-h} \right)^n}{\left[ n \right]!}  \left( E_- \otimes E_+
\right)^n e^{\frac{1}{4}h n(n-1)}, \qquad \left[ n \right] = \frac{e^{\frac{%
nh}{2}}-e^{\frac{-nh}{2}}}{e^{\frac{h}{2}}-e^{\frac{-h}{2}}}.
\end{equation}
The same is true for the hybrid algebra $U_{q{\cal J}} (sl(2))(h,\xi)$.
According to the general properties of  twisted quasitriangular algebras
(see eq. (\ref{RF})) $U_{q{\cal J}} (sl(2))(h,\xi)$ has the following ${\cal %
R}$-matrix,
\begin{equation}
\label{uni-rqj} {\cal R}_{q{\cal J}} = e^{\omega \otimes H} e^{-hH \otimes
H} {\cal R}_q e^{hH \otimes H} e^{- H \otimes \omega}.
\end{equation}
%%%%%%%%%%%%%%%%%%%%%%%%%%
In a case of the smooth set of quantized algebras the classical limit
depends on how we fix the linear subvariety that describes the deformation
quantization. If we want to disclose the hybrid properties of the set $%
\left\{ U_{q{\cal J}} (sl(2))(h,\xi) \right\} $ we are to find a smooth
curve intermediate between the standard deformation subvariety $\left\{ U_{q%
{\cal J}} (sl(2))(h,0) \, | \, h \geq 0 \right\} $ and the pure twist
subvariety $\left\{ U_{q{\cal J}} (sl(2))(0,\xi) \, | \, \xi \geq 0 \right\}
$. Obviously, it is sufficient to consider a linear subvariety $\left\{ U_{q%
{\cal J}} (sl(2))(\zeta \xi,\xi)\,| \,\xi \geq 0, \, \zeta > 0 \right\} $
where we had put $h= \zeta \xi$. In the corresponding set $\left\{{\cal R}_{q%
{\cal J}} \, | \,h= \zeta \xi \right\}$ of hybrid ${\cal R}$-matrices (\ref
{uni-rqj}) we let $\xi$ to be in the neighborhood of zero and extract the
classical $r$--matrix
\begin{equation}
\label{r-qj} r_{q{\cal J}} = E_+ \wedge H + \zeta \left( H \otimes H + E_-
\otimes E_+ \right) .
\end{equation}
This expression is the well known hybrid solution \cite{GGS} of the
classical Yang--Baxter equation. 
%%%%%%%%%%%%%%%%%%%%%%%%%%
%%%%%%%%%%%%%%%%%%%%%%%%%%

\subsection{Hybrid quantum affine algebra $U_{q{\cal J}}(\widehat{sl(2))}$}

The explicit construction of Jordanian twist \cite{OGI}, extended Jordanian
twist \cite{KLM} and chains of twists \cite{KLO} provided the possibility to
obtain the mixed quantizations for current algebras -- the twisted Yangians
\cite{KS,L-DUB,KLS-PR}. Analogously, with the help of the $q$--Jordanian
twist ${\cal F}_{q{\cal J}}(h,\xi)$ we can obtain the hybrid quantizations
for Kac--Moody algebras.

Let us consider, for example, the quantum affine algebra $U_{q}(\widehat{sl(2))%
}$ \cite{DRI2,JIM} defined as a deformed infinite dimensional Lie algebra
with the Cartan matrix
$$
A=(a_{ij})= [\left( \lambda_i , \lambda_j \right)] = \left(
\begin{array}{rr}
2 & -2 \\
-2 & 2 \\
\end{array}
\right), \qquad i,j = 0,1;
$$
the generators $H_i$, $E_{\pm\lambda_i}$, $D$ and the relations
\begin{equation}
\label{aff-r1}
\begin{array}{l}
\left[ H_i,E_{\pm\lambda_j} \right] = \pm
\frac{1}{2} a_{ij} E_{\pm\lambda_j}, \\ [0.2cm] \left[
E_{\lambda_i},E_{-\lambda_j} \right] = \delta_{ij}
\frac{e^{hH_i}-e^{-hH_i}}{1 - e^{-h}}, \\ [0.2cm] \left[ D,E_{\pm\lambda_i}
\right] = \pm \delta_{i0}E_{\pm\lambda_i}, \qquad i,j = 0,1; \\
[0.2cm] \left[ H_i,H_j \right] = \left[ H_i,D \right] = 0, \\
[0.2cm] \left( {\rm ad}_q E_{\pm\lambda_i} \right)^{1 - a_{ij}}\circ
E_{\pm\lambda_j} = 0 , \qquad i \neq j.
\end{array}
\end{equation}
Here ${\rm ad}_q$ is the $q$--adjoint operator
$$
{\rm ad}_q E_{\lambda_i}\circ E_{\lambda_j} = E_{\lambda_i} E_{\lambda_j} -
e^{h \left( \lambda_i , \lambda_j \right)} E_{\lambda_j} E_{\lambda_i}.
$$
We shall put $q= e^{\frac{1}{2}h}$ and introduce the rescaled generators
\begin{equation}
\label{e-gen} {\bf e}_{\pm\lambda_i} = e^{-\frac{1}{4}h}E_{\lambda_i}.
\end{equation}
Let $\delta = \lambda_0 + \lambda_1$ be the minimal imaginary root of $
\widehat{ sl(2)}$. Then, the so called normal ordering \cite{KT1} in the
system of positive roots $\Lambda_+$ is fixed as follows
\begin{equation}
\label{nor-ord} \lambda_0,\lambda_0 +\delta, \dots , \lambda_0 + n\delta,
\dots , \dots, \delta ,\ 2 \delta , \dots, \dots, \lambda_1 + (l+1)\delta ,\
\lambda_1 + l \delta, \dots , \lambda_1.
\end{equation}
According to this ordering the generators for composite roots are obtained
as follows
$$
\begin{array}{lr}
{\bf e}^{\prime}_{\delta} = \left[ 2 \right]^{-1} \left[ {\bf e}%
_{\lambda_0}, {\bf e}_{\lambda_1} \right]_q, \qquad {\bf e}%
^{\prime}_{\lambda_0 + n\delta} = (-1)^n \left( {\rm ad} {\bf e}%
^{\prime}_{\delta}\right)^n \circ {\bf e}_{\lambda_0}, & \\[0.2cm]
{\bf e}^{\prime}_{\lambda_1 + n\delta} = \left( {\rm ad} {\bf e}%
^{\prime}_{\delta}\right)^n \circ {\bf e}_{\lambda_1}, \quad {\bf e}%
^{\prime}_{n\delta} = \left[ 2 \right]^{-1} \left[ {\bf e}_{\lambda_0 +
(n-1)\delta}, {\bf e}_{\lambda_1} \right]_q. \\  &
\end{array}
$$
(The $q$-numbers above are the same as in (\ref{uni-r}).)  Finally, the
generators ${\bf e}_{n\delta}$ are defined by means of the Schur
polynomials:
$$
{\bf e}^{\prime}_{n \delta} = \sum_{p_1+2p_2+ \dots +np_n=n} \frac{%
(q^2-q^{-2})^{\sum p_i-1}}{p_1! \dots p_n!} {\bf e}_{\delta}^{p_1} {\bf e}%
_{2\delta}^{p_2} \dots {\bf e}_{n\delta}^{p_n}.
$$
The generators for the negative roots are defined with the help of the
involution
$$
\left( H_i \right) ^* = - H_i, \qquad \left({\bf e}_{\pm\lambda_i} \right)
^* = {\bf e}_{\mp\lambda_i}, \qquad h^* = -h.
$$
In term of these generators the universal ${\cal {R}}$-matrix of $U_q (
\widehat{sl(2)})$ has the form \cite{KT3}
\begin{equation}
\label{rmat1}
\begin{array}{c}
{\cal R}^{DJ} = \left( \prod_{n \leq 0}^{\longrightarrow} \exp_{q} \left( (q
- q^{-1}) {\bf e}_{\alpha + n \delta} \otimes {\bf e}_{- \alpha - n \delta}
\right) \right) \cdot \exp \left( \sum_{n>0} \frac{n({\bf e}_{n \delta}
\otimes {\bf e}_{-n \delta})}{ q^{2n} - q^{-2n}} \right) \\ [0.2cm] \cdot
\left( \prod_{n \leq 0}^{\longleftarrow} \exp_{q} ((q - q^{-1}) {\bf e}%
_{\beta + n \delta} \otimes {\bf e}_{- \beta - n \delta} ) \right) \cdot
{\cal K},
\end{array}
\end{equation}
where ${\cal K}$ stands for
$$
{\cal K} = \exp \left( \sum_{i,j} 2hd_{ij} H_i \otimes H_j \right),
$$
$d$ is the inverse of the extended (nondegenerate) Cartan matrix $\widetilde{%
a} $ \cite{Kac} and the $q$--exponent is defined as the series
$$
\exp_q \equiv \sum \frac{x^n}{(n)_{q^{-2}}}, \qquad (n)_{q^{-2}} = \frac{%
q^{-2n}-1}{q^{-2}-1}.
$$
Note that the order of $q$-exponents in (\ref{rmat1}) is direct in the first
product ($\rightarrow$) and inverse in the second one ($\leftarrow$).

Any quantum Borel subalgebra $U_q (B) \in U_q (\widehat{sl(2)}) $ can be
used as a carrier algebra to perform the quantum Jordanian twisting. If we
have to consider representations of the corresponding quasitriangular
quantum algebras the simplest choice is to take the Hopf subalgebra
generated by $H_0$ and $E_{\lambda_0}$. The twist deformation is performed
by the element (see (\ref{qr-twist}))
\begin{equation}
\label{qr-twists} {\cal F}_{q{\cal J}}(h,\xi) = e^{H_0 \otimes \omega_0}
e^{-hH_0 \otimes H_0}, \qquad \omega_0 = \ln \left( \left( \xi E_{\lambda_0}
+ 1 \right) e^{hH_0} \right)
\end{equation}
and produces the Jordanian quantum affine algebra $U_{q{\cal J}}(\widehat{%
sl(2))}$. It has the commutators defined by (\ref{aff-r1}) and the deformed
coproducts:
\begin{equation}
\begin{array}{lcl}
\Delta_{q{\cal J}} (H_i) & = & H_i \otimes 1 + e^{
{\rm ad} (H_0 \otimes w)} \circ (1 \otimes H_i); \\ [0.2cm] \Delta_{q{\cal J}%
} (D) & = & D \otimes 1 + e^{
{\rm ad} (H_0 \otimes w)} \circ (1 \otimes D); \\ [0.2cm] \Delta_{q{\cal J}}
(E_{ \lambda_0}) & = & \frac{1}{\xi} (e^w \otimes e^w - e^{{\rm ad} (H_0
\otimes w)} \circ (e^{-hH_0} \otimes e^{-hH_0}) -1 \otimes 1); \\ [0.2cm]
\Delta_{q{\cal J}} (E_{- \lambda_0}) & = & E_{- \lambda_0} \otimes e^{-w} +
e^{
{\rm ad} (H_0 \otimes w)} \circ (1 \otimes E_{- \lambda_0}); \\ [0.2cm]
\Delta_{q{\cal J}} (E_{ \lambda_1}) & = & (E_{ \lambda_1} \otimes e^{-w})
\cdot ( e^{
{\rm ad} (H_0 \otimes w)} \circ (1 \otimes e^{hH_0})) + e^{h(H_1 + H_0)}
\otimes E_{ \lambda_1}; \\ [0.2cm] \Delta_{q{\cal J}} (E_{- \lambda_1}) & =
& (E_{- \lambda_1} \otimes e^{+w}) \cdot (e^{
{\rm ad} (H_0 \otimes w)} \circ (1 \otimes e^{-h(H_1 + H_0)})) \\ [0.1cm] &
& \qquad\qquad + e^{-hH_0} \otimes E_{- \lambda_1}.
\end{array}
\end{equation}

The twisted (hybrid) universal ${\cal R}$--matrix for $U_{q{\cal J}}(
\widehat{sl(2))}$ has the following form
\begin{equation}
\label{rmat2}
\begin{array}{ll}
{\cal R}^{DJ}_{q{\cal J}} = & e^{\omega \otimes H} e^{-hH \otimes H} \cdot
\left( \prod_{n \leq 0}^{\longrightarrow} \exp_{q} \left( (q - q^{-1})
{\bf e}_{\alpha + n \delta} \otimes {\bf e}_{- \alpha - n \delta} \right)
\right) \\ [0.2cm] & \quad \cdot \exp \left( \sum_{n>0}
\frac{n({\bf e}_{n \delta} \otimes {\bf e}_{-n \delta})}{ q^{2n} - q^{-2n}}
\right) \cdot \left( \prod_{n \leq 0}^{\longleftarrow} \exp_{q} ((q -
q^{-1}) {\bf e}_{\beta + n \delta} \otimes {\bf e}_{- \beta - n \delta} )
\right) \\ [0.2cm] & \quad \cdot {\cal K} \cdot e^{hH \otimes H} e^{- H
\otimes \omega}.
\end{array}
\end{equation}
It satisfies the parametric QYBE
$$
{\cal R}_{12} (z_1/z_2) {\cal R}_{13} (z_1/z_3) {\cal R}_{23} (z_2/z_3) =
{\cal R}_{23} (z_2/z_3) {\cal R}_{13} (z_1/z_3) {\cal R}_{12} (z_1/z_2).
$$
In the fundamental representation of $sl(2)$ we get the hybrid matrix
solution:
\begin{equation}
\label{drm}
\begin{array}{c}
d \left( {\cal R}_{q{\cal J}}^{ DJ} \right) = \frac{1-z}{g^{\frac{3}{2}}
(1-zq^{-2})} \exp \left( \sum_{n>0} \frac{z^n}{n} \frac{q^n - q^{-n}}{q^n +
q^{-n}} \right) \cdot \left(
\begin{array}{cccc}
a_1 & sq & -s & s^2 \\
0 & q & a_2 & s \\
0 & za_2 & q & -sq \\
0 & 0 & 0 & a_1 ,
\end{array}
\right)
\end{array}
\end{equation}
where
$$
a_1 = \frac{q^2 -z}{1-z} , \qquad a_2 = \frac{q^2 -1}{1-z} , \qquad s =
\frac{\xi}{1+q}.
$$

The expression for the universal ${\cal R}$-matrix (\ref{rmat2}) as well as
for its defining representation (\ref{drm}) describes a smooth variety of
solutions of QYBE. This is the two-dimensional variety with the coordinates $%
h$ and $\xi$ and with the spectral parameter $z$:
\begin{equation}
\begin{array}{ccccc}
&  & {\cal R}_{q{\cal J}}^{DJ} (h,\xi) &  &  \\
& \rule{25mm}{0mm} {\scriptsize h \rightarrow 0} \swarrow &  & \searrow
{\scriptsize \xi \rightarrow 0}\rule{25mm}{0mm} &  \\
& {\cal R}_{{\cal J}} = {\cal R}_{q{\cal J}}^{DJ} (0,\xi) &  & {\cal R}_{q%
{\cal J}}^{DJ} (h,0) =  {\cal R}^{DJ} &
\end{array}
\end{equation}
When $\xi$ goes to zero we return to the initial quantum affine algebra $%
U_{q}(\widehat{sl(2)})$ and the corresponding ${\cal R}$-matrix (\ref{rmat1}%
). In the limit $h \rightarrow 0 $ we get the nonstandard quantization $U_{%
{\cal J}}(\widehat{sl(2)})$ of the affine algebra $U (\widehat{sl(2)})$
performed by the Jordanian twist
$$
{\cal F_{J}} = e^{H_0 \otimes \sigma_0} ,
$$
with $\sigma_0 = \ln \left( 1 + E_{\lambda_0} \right) $. The ${\cal R}$%
-matrix in this limit case becomes the ordinary Jordanian.

%%%%%%%%%%%%%%%%%%%%%%%%%%%%%%%%%%%%%%%%%%%%%%%%%%%%%%%%%%%%%%%%%%%%%%%%%
%%%%%%%%%%%%%%%%%%%%%%%%%%%%%%%%%%%%%%%%%%%%%%%%%%%%%%%%%%%%%%%%%%%%%%%%%

\sect{Conclusions}

\label{conclusions}

We have demonstrated that there exists the hybrid quantization $U_{q{\cal J}%
}(sl(2))$ with well defined natural limits with respect to the two
deformation parameters $h$ and $\xi $. Moreover, each quantum algebra $U_{q%
{\cal J}}(sl(2))$ can be considered together with the twisting element $%
{\cal F}_{q{\cal J}}$ that connects it with the corresponding standard
quantization $U_q(sl(2))$. Both natural limits exist also for the triples $%
(U_{q{\cal J}}(sl(2)),\ {\cal F}_{q{\cal J}},\ {\cal R}_{q{\cal J}}^{DJ})$.
Such limit behaviour illustrates the difference between previously obtained
combined quantizations \cite{BHP,CHA,ABD,KTS} and the deformation $U_{q{\cal %
J}}$ produced by the quantum Jordanian twist. Contrary to the cases of $(q,\
\xi )$--deformation by Ballesteros {\sl et al} \cite{BHP} (BHP--deformation)
and the constructions proposed by Abdesselam {\sl et al} \cite{CHA,ABD} and
Stolin \cite{KTS} in the triple $(U_{q{\cal J}}(sl(2)),\ {\cal F}_{q{\cal J}%
},\ {\cal R}_{q{\cal J}}^{DJ})$ there are no singularities when $q\to 1$.
This means that the sets $\{U_{q{\cal J}}(sl(2))\}$ and $\{(U_{q{\cal J}%
})(sl(2))_{{\cal J}}\}$, as well as the BHP--deformation, refer to different
``sheets''.

In a review by Aneva {\sl et al} \cite{CHA2} the BHP--deformation was
considered as being not an authentic hybrid quantization. The reason was
that in the generic points the standard and the BHP--deformations are
equivalent (as it was proved in \cite{CHA2}). From our point of view the
requirement of nonequivalence to $U_q$ is too strong in the context of
hybrid quantizations. The main criterion here must be the possibility for a
Hopf algebra to be submerged in a two-dimensional smooth variety whose
boundaries are $U_{q}$ and $U_{{\cal J}}$. As it was shown in \cite{KL} this
implies that it is a quantization of a hybrid classical $r$--matrix (\ref
{r-qj}).

Notice that the Jordanian quantum affine algebra $U_{q{\cal J}}(\widehat{%
sl(2))}$ constructed in Section \ref{jordanian
deformations} is an example
of {\sl twisted quantum nontwisted affine algebras}, $U_{q{\cal J}}(\widehat{%
sl(2))} = U_{q{\cal J}}(A_1^{(1)})$. The word ``twisted" in the term
``twisted affine algebra" (introduced by Kac \cite{Kac}) has the meaning
different from that of the Drinfeld's deformation procedure \cite{DRI}. This
is why the term ``Jordanian" is preferable here.

The quantum Jordanian twist ${\cal F}_{q{\cal J}}$ can be applied to any
Hopf algebra containing the quantum Borel subalgebra $U_q(B)$. In particular
it can be used to produce hybrid deformations of the twisted affine algebras
(the Kac--Moody algebras listed in the tables Aff2 and Aff3 in \cite{Kac}).

%%%%%%%%%%%%%%%%%%%%%%%%%%%%%%%%%%%%%%%%%%%%%%%%%%%%%%%%%%%%%%%%%%%%%%

\section*{Acknowledgments}

V. L. would like to thank the Vicerrectorado de Investigaci\'on de la
Universidad de Va\-lla\-dolid for supporting his stay. This work has been
partially supported by DGES of the Ministerio de Educaci\'on y Cultura of
Espa\~na under Project PB98-0360, the Junta de Castilla y Le\'on (Espa\~na)
and the Russian Foundation for Fundamental Research under the grant
00-01-00500.

%%%%%%%%%%%%%%%%%%%%%%%%%%

%%%%%%%%%%%%%%%%%%%%%%%%%%

%%%%%%%%%%%%%%%%%%%%%%%%%

\end{document}